\documentclass[12pt,oneside]{amsart}
\usepackage[a4paper,scale=.76]{geometry}
\usepackage{natbib}
\usepackage[colorlinks,citecolor=blue,urlcolor=blue]{hyperref}
\usepackage{graphicx}
\usepackage{amsfonts,amsmath,amssymb,bbm}
\def\proof{\noindent{\bf Proof:}\hskip10pt}        
\def\QED{\hfill $\Box$}

\font\tenmath=msbm10 scaled 1200
\font\sevenmath=msbm7 scaled 1200
\font\Fivemath=msbm5 scaled 1200
\newfam\mathfam \textfont\mathfam=\tenmath
\scriptfont\mathfam=\sevenmath \scriptscriptfont\mathfam=\Fivemath

\def \\ { \cr }
\def\R{\mathbb{R}}
\def \1{1 \mkern -6mu 1} 
\def\N{\mathbb{N}}
\def\E{\mathbb{E}}
\def\P{\mathbb{P}}
 
\def\Z{\mathbb{Z}}
\def\R{\mathbb{R}}

\def \e{{\rm e}}
\def \d{{\rm d}}
\def\x{{\tt x}}

\newtheorem{theorem}{Theorem}

\newtheorem{lemma}{Lemma}
\newtheorem{corollary}{Corollary}
\begin{document}

\title[Percolation on scale-free random trees]{Supercritical percolation \\ on large scale-free random trees}

\author{Jean Bertoin}
\address[Jean Bertoin]{Institut f\"ur Mathematik, 
Universit\"at Z\"urich, 
Winterthurerstrasse 190, 
CH-8057 Z\"urich} 
\author{Ger\'onimo Uribe Bravo}
\address[Ger\'onimo Uribe Bravo]{Instituto de Matem\'aticas, Ciudad Universitaria
Coyoac\'an, 04510. M\'exico, D. F. }

\date{\today}
\begin{abstract}  
We consider Bernoulli bond percolation on a large scale-free tree in the supercritical regime, meaning informally  that  there exists a giant cluster with high probability.
We obtain a weak limit theorem for the sizes of the next largest clusters, extending a recent result in \cite{BERTOIN:2012:HAL-00636264:1} for large random recursive trees. The approach relies on 
 the analysis of the asymptotic behavior of branching processes subject to rare neutral mutations, which may be of independent interest.

\end{abstract}
\maketitle
\begin{section}{Introduction and statement of the main result}
Generalizing the procedure to grow graphs and trees of \cite{MR2091634}  (see also \cite{MR930497} for an earlier version),  \cite{PhysRevLett.85.4633} and \cite{MR1909153}
 grow a so-called scale-free random tree on a set of ordered vertices, say  $\{0,\ldots,n\}$, using an algorithm with preferential attachment that we now recall.   
 Fix a parameter
$\beta\in(-1, \infty)$, and start  for $n=1$ from  
the unique tree $T_1$ on $\{0,1\}$ which has a single edge connecting $0$ and $1$. Then suppose that $T_n$ has been constructed for some $n\geq 1$, and for every $i\in\{0,\ldots,n\}$, denote by $d_n(i)$ the degree of the vertex $i$ in $T_n$. Conditionally given $T_n$, the tree $T_{n+1}$ is derived from $T_n$ by incorporating  the new vertex $n+1$ and creating an edge  between $n+1$ and a vertex $v_n\in\{0,\ldots,n\}$ chosen at random according to the law
$$\P(v_n=i\mid T_n)= \frac{d_n(i)+\beta}{2n+\beta(n+1)}\,, \qquad i\in\{0,\ldots,n\}.$$
That the preceding indeed defines a probability on $\{0,\ldots,n\}$ for $\beta<\infty$  is seen from the fact that, since $T_n$ is a tree with $n+1$ vertices, there is the identity $\sum_{i=0}^nd_n(i)=2n$. 

There has been a significant interest in the last decade in Bernoulli bond-percolation on large scale-free graphs; see in particular \cite{MR2076725} and \cite{MR2174664}.
In the simpler case of trees, this means that having constructed $T_n$ for some $n\gg 1$, and for a given parameter $p(n)\in(0,1)$, we keep each edge with probability $p(n)$ and remove it with probability $1-p(n)$, independently of the other edges. This disconnects $T_n$ into a family of clusters, and the purpose of this work is to study the asymptotic behavior in distribution of the sizes of the largest clusters as $n\to\infty$, for a particular regime of the sequence $p(n)$. 
Specifically, let us write
$$C^{(p)}_{0,n}\geq C^{(p)}_{1,n}\geq \ldots$$
for the ordered sequence of the sizes of the clusters. 

In the boundary case $\beta\to \infty$ where $v_n$ becomes uniformly distributed on $\{0,\ldots,n\}$, the algorithm yields a so-called uniform recursive tree (see for instance \cite{MR1445048,MR2484382}). 
It has then been observed recently by \cite{BERTOIN:2012:HAL-00636264:1}, that choosing the percolation parameter so that 
\begin{equation}\label{Eq1}
1-p(n)\sim  \frac{c}{\ln n}
\end{equation}
where $c>0$ is fixed, corresponds precisely to the supercritical regime, in the sense that  both the largest percolation cluster  on a random recursive tree of size $n\gg 1$ and its complement, have a size of order $n$ with high probability. Specifically, the largest cluster has a size close to $\e^{-c}n$ whereas the next largest clusters have size of order $n/\ln n$ only and are approximately distributed according to some Poisson random measure with intensity $c\e^{-c} x^{-2} \d x$.  

The main purpose of this work is to show that a similar result holds more generally for large scale-free random  trees. 

\begin{theorem}\label{T1} Set $\alpha = (1+\beta)/(2+\beta)$, and assume that the percolation parameter $p(n)$ fulfills \eqref{Eq1}. Then 
$$\lim_{n\to \infty}n^{-1}C^{(p)}_{0,n}=\e^{-\alpha c }\qquad \hbox{in probability,}$$
and for every fixed integer $j$, 
$$\left(\frac{\ln n}{n}C^{(p)}_{1,n}, \ldots, \frac{\ln n}{n}C^{(p)}_{j,n}\right)$$
converges in distribution towards
$$(\x_1, \ldots, \x_j)$$
where $\x_1>\x_2>\ldots$ denotes the sequence of the atoms of a Poisson random measure
on $(0,\infty)$ with intensity 
$$c \e^{-\alpha c } x^{-2}\d x\,.$$
Equivalently,  $1/\x_1, 1/\x_2-1/\x_1, \ldots, 1/\x_j-1/\x_{j-1}$ are i.i.d. exponential variables with parameter $c \e^{-\alpha c }$. 

\end{theorem}

It is remarkable that the intensity measure in the statement only depends on the parameter $\beta$ through the constant factor $\e^{-\alpha c }$. 
It should also be noted that  the map $\beta\mapsto \alpha(\beta)=(1+\beta)/(2+\beta)$ increases, and we then see from Theorem \ref{T1} that  for the same value of the percolation parameter $p(n)$ and $n\gg 1$, this intensity decreases with the parameter $\beta$. This can be explained informally by the fact that  when  the parameter $\beta$ is larger, the algorithm with preferential attachment produces random trees which are less tufty  and thus more affected by percolation.

The approach used in \cite{BERTOIN:2012:HAL-00636264:1} for recursive trees relies crucially on special properties of the latter, and more specifically on a remarkable coupling due to  \cite{MR2407414} connecting the Meir and Moon algorithm for the isolation of the root with a certain random walk in the domain of attraction of the completely  asymmetric Cauchy process. This fails for scale-free trees, and we thus have to use here a fairly different route. 

It is well-known that growing random scale-free trees bears close relations to Yule processes. We shall incorporate an independent Bernoulli percolation to the algorithm with preferential attachment and interpret this in terms of neutral mutations which are superposed to the structure of the branching process.
This leads us to investigate in Section \ref{branchingSection}  the asymptotic behavior of a system of branching processes with rare neutral mutations up to a large random time, in certain regimes when the small mutation parameter is related to the size of the total population.  We then specify in Section \ref{percolationSection} those results to Yule processes, make the link with percolation on scale-free trees and prove Theorem \ref{T1}.
\end{section}

\begin{section}{Branching processes with rare neutral mutations}
\label{branchingSection}
Thus main purpose of this section is to establish some general results about the long time behavior of a system of branching processes with rare neutral mutations in a certain specific regime. The system is presented in the first sub-section, and then asymptotic results are established in the second. 

\subsection{Description of the system of branching processes with mutations}

We start by considering a pure birth branching process $Z=(Z(t): t\geq 0)$ in continuous space, with unit birth rate per unit population size and reproduction law $\nu$, where $\nu$ denotes a probability measure on $(0,\infty)$. This means that $Z$ is a non-decreasing  Markovian jump process such that 
when $Z(0)=z>0$, its first jump occurs after an exponential time with parameter $z$, and the jump size has  law $\nu$.

We assume that the second moment of $\nu$ is finite, which is more than sufficient to ensure that $Z$ never explodes a.s. Recall that $\beta >-1$ is some fixed parameter. We further suppose that $\nu((0,1+\beta])=0$ (the role of this assumption shall be plain latter on), so that when a birth even occurs, the population always increases by an amount at least  $1+\beta$.
We shall be mainly interested  in a class of 
population systems which arise by incorporating neutral mutations to the preceding branching process. 
It may be useful to think of Kimura's infinite site model, in which a genetic type consists of an infinite sequence of letters and  each mutation affects  a different locus. In particular, one can reconstruct the genealogy of the types by comparison with the ancestral type; see e.g. Section 2 in  \cite{MR2603059}
for a closely related setting.

More precisely, let  ${\mathbf U}=\bigcup_{n\geq 0}\N^n$ denote the Ulam tree, with the convention that $\N^0=\{\varnothing\}$.
That is, each element $u\in {\mathbf U}$ is a finite sequence $u=(u_1, \ldots, u_n)$ of positive integers, whose length $|u|=n$ corresponds to the height of $u$ in $ {\mathbf U}$, 
and the empty sequence $\varnothing$ serves as the root of ${\mathbf U}$. Each vertex $u\in {\mathbf U}$ corresponds to a genetic type; in particular we view $\varnothing$ as the ancestral type, and for every  $u=(u_1, \ldots, u_n)\in {\mathbf U}$ and $j\in\N$, the $j$-th child of $u$, $uj= (u_1, \ldots, u_n, j)$, represents  the new genetic type which appears at the instant when the $j$-th mutation occurs in the sub-population with type $u$.

The state of the population system at a given time is given by
a collection of nonnegative real numbers $(z_u: u\in {\mathbf U})$, where $z_u$ is the current size of the sub-population with  type $u$. The evolution of the system is thus described  by a process
${\bf Z}=({\bf Z}(t), t\geq 0)$, where for each $t\geq 0$, ${\bf Z}(t)= (Z_u(t): u\in {\mathbf U})$ is a collection of nonnegative variables indexed by Ulam's tree. 
At the initial time, all the $Z_u(0)$ are taken to be equal to zero, except $Z_{\varnothing}(0)$ which is the size of the ancestral population. 

We then describe the random evolution of the system ${\bf Z}$, which depends on  a parameter $p\in[0,1]$.
Recall that the reproduction law $\nu$ assigns no mass to $(0,1+\beta]$, so we may consider a positive random variable $\xi$ such that $\xi + 1+\beta$ has the law $\nu$.
We imagine that mutations occur at rate $1-p$ per unit population size,  always produce a single mutant population of fixed size $1+\beta$, and are neutral, in the sense that they do not affect the reproduction law. 
In particular,  the different populations present in the systems (i.e. with strictly positive sizes) evolve independently one of the other and according to the same random dynamics. 
For each  sub-population, say with size $z>0$, we   
introduce an independent copy $\xi'$ of $\xi$, a Bernoulli variable $\epsilon_p$  with parameter $p$, and an exponentially distributed variable  $\zeta_z$ with parameter $z$. We assume that these three variables are independent, and also independent of the other variables associated to the other sub-populations. 
The  time $\zeta_z$ corresponds to the first birth event in that sub-population. The total size of the children born at this birth event is $\xi'+1+\beta$; the 
 variable $\epsilon_p$ specifies whether a mutation occurs. Specifically, the size of clone children is $\xi'+\epsilon_p (1+\beta)$
 and the size of mutant children is  $(1-\epsilon_p)(1+\beta)$. So mutation occurs if and only if $\epsilon_p = 0$, an event which has probability $1-p$. 
 
 In order to underline the role of the rate of mutation, we henceforth write 
 $${\bf Z}^{(p)}=(Z^{(p)}_u(t): t\geq 0, u\in {\mathbf U})$$
 instead of ${\bf Z}$.  It should be obvious however that, no matter what $p$ is, the process of the total size of the population
 $$Z(t) = \sum_{u\in{\mathbf U} } Z^{(p)}_u(t)\,,\qquad t\geq 0\,,$$
 is distributed as the branching process described at the beginning of this section.

Clearly, the process of the size of the sub-population with the ancestral type $(Z^{(p)}_{\varnothing}(t): t\geq 0)$ is  a continuous time branching process in continuous space with reproduction law given by the distribution of 
$ \epsilon_p(1+\beta)+ \xi$. More generally, 
if for $u\in{\mathbf U}$, we write 
$${b}^{(p)}_u=\inf\{t\geq 0: Z^{(p)}_u(t)>0\}\,,$$
for the birth time of the sub-population with type $u$, 
then each process 
$$(Z^{(p)}_u(t+{b}^{(p)}_u), t\geq 0)$$
is a branching processes with the same reproduction law as $Z^{(p)}_{\varnothing}$ and starting from $1+\beta$ for $u\neq \varnothing$. 
Further, it should be intuitively clear (although this shall not be needed here) that the processes 
$Z^{(p)}_u({b}^{(p)}_u+ \cdot)$ for $u\in{\mathbf U}$ are independent. Focussing  on types of the first generation $\N^1$, i.e. bearing a single mutation, we point at a useful independence property involving the birth times: 

\begin{lemma} \label{Le2} The processes  $(Z^{(p)}_i({b}^{(p)}_i+t): t\geq 0)$ for $i\geq 1$ form a sequence  of i.i.d. branching processes with reproduction distributed according to $\xi+\epsilon_p (1+\beta)$,
and starting point $1+\beta$. Further, this sequence is independent of that of
the birth-times  $({b}^{(p)}_i)_{i\geq 1}$ and of the process   $Z^{(p)}_{\varnothing}$ of the sub-population with the ancestral type. 
\end{lemma}

\noindent{\bf Remark.} It is crucial in this statement to focus on sub-populations bearing the same number of mutations; for instance  the independence property of the birth-times would  fail if we considered the whole the family of processes $Z^{(p)}_u({b}^{(p)}_u+\cdot)$ for $u\in{\mathbf U}\backslash\{\varnothing\}$.

\proof 
Let $(X,M)$ be a continuous-time Markov chain with values in $\R_+\times \Z_+$ with two types of transition:
\begin{align*}
(x,m) \mapsto\ (x+\d x,m)&\text{ at rate }x p \nu(\d x)\\
(x,m) \mapsto\ (x+\d x,m+1)&\text{ at rate }x(1-p)\nu( 1+\beta+\d x).
\end{align*}
In particular, $X$ is a branching process distributed as $Z^{(p)}_{\varnothing}$ and we can interpret $M$ as the process of the number of mutation events which occur within the sub-population with the ancestral type.

Let $\gamma_1<\gamma_2<\cdots$ denote the sequence of  jump times of $M$ and set $\gamma_0=0$. Independently on $(X,M)$, let $(X_i, i\in\N)$ be a sequence of  i.i.d. branching processes with the same law as $Z^{(p)}_{\varnothing}$ but with starting value $1+\beta$. 
We then form the process 
$${\bf X}(t)=\left(X(t), {\bf 1}_{t\geq \gamma_1}X_1(t-\gamma_1), {\bf 1}_{t\geq \gamma_2}X_2(t-\gamma_2), \ldots \right)\,,\qquad t\geq 0\,.$$
The analysis of jump times and positions then readily shows that ${\bf X}$ is Markovian and has the same law as $(Z^{(p)}_{\varnothing}, Z^{(p)}_1, Z^{(p)}_2, \ldots)$.
\QED

\subsection{Asymptotics for rare mutations}
Recall the assumption  that the reproduction law $\nu$ of the branching process $Z$ has a finite second moment and write
$$m_1=\int x \nu(\d x) \quad \hbox{and}\quad m_2=\int x^2\nu(\d x) \,.$$ 
It is well-known that 
$$W(t):=\e^{-m_1t}Z(t)\,,\qquad t\geq 0$$
is then a nonnegative square-integrable martingale, and we write $W(\infty)$ for its terminal value. Furthermore $W(\infty)>0$ a.s. since $Z$ cannot become extinct (cf. Theorem 2 p. 112 in Chapter III of \cite{MR2047480} for the general assertion and Example 5.4.3 p. 253 of \cite{MR2722836} just for the finite variance case).

It is easily checked that the speed of convergence of the martingale $W$  is exponential. 
Specifically, if we write $\P_z$ for the distribution of the branching process $Z$ started from $z>0$, then the following general bound
holds\footnote{The assumption that $\nu$ assigns no mass to $(0,1+\beta]$ plays no role here, and Lemma \ref{L1} holds when this assumption is dropped.}.

\begin{lemma}\label{L1} For every $t\geq 0$, there is the upper-bound
$$\E_z\left(\sup_{s\geq t}|W(s)-W(\infty)|^2 \right) \leq 10z \frac{m_2}{m_1} \e^{-m_1t}\,.
$$
As a consequence, we have 
$$\E_z\left(\sup_{s\geq 0}\e^{m_1 2s/3}|W(s)-W(\infty)|^2 \right) \leq  \frac{10z m_2 \e^{2m_1/3}}{m_1
(1-\e^{-m_1/6})^2} 
\,.$$
\end{lemma} 

\proof By Doob's inequality and basic properties of square integrable martingales, we have
$$\E_z\left(\sup_{s\geq t}|W(s)-W(\infty)|^2 \right) \leq 10 \E_z([W]_{\infty}-[W]_t)$$
where
$$[W]_t=\sum_{0\leq s \leq t} |\Delta W(s)|^2 =\sum_{0\leq s \leq t} \e^{-2m_1s} |Z(s)-Z(s-)|^2\,.$$
A straightforward calculation shows that the compensator of jump process $[W]$ is
$$\langle W\rangle_t= m_2 \int_{0}^t\e^{-2m_1s}Z(s) \d s\,,$$
that is $[W]_t-\langle W\rangle_t$ is a local martingale. 
Finally observe that 
$$\E_z(\e^{-2m_1s}Z(s))= \e^{-m_1s}\E_z(W(s))=z\e^{-m_1s}\,,$$
so 
$$\E_z(\langle W\rangle_{\infty}-\langle W\rangle_t)= z \frac{m_2}{ m_1} \e^{-m_1t}\,.$$ 
This enables us to assert that 
$$ \E_z([W]_{\infty}-[W]_t) =  \E_z(\langle W\rangle _{\infty}-\langle W\rangle _t)\,,$$
and our first statement follows. 

Turning our attention to the second inequality, we  write for every integer $n\geq 0$
$$\sup_{n\leq s< n+1}\e^{m_1 s/3}|W(s)-W(\infty)| \leq \e^{m_1 (n+1)/3}\sup_{n\leq s< n+1} |W(s)-W(\infty)|\,.$$
It follows from the first part that the $L^2$-norm of the right hand side can be bounded from above by 
$$\e^{m_1 (n+1)/3}\|\sup_{s\geq n} |W(s)-W(\infty)|\|_2\leq \sqrt{10z \frac{m_2}{m_1}} \e^{-m_1(n-2)/6}\,,
$$
so taking the sum over $n$ and applying Minkowski's inequality yields the stated bound. 
\QED

The main purpose of this section is to specify  the joint asymptotic behaviors of  the branching processes $Z^{(p)}_{\varnothing}$ and $Z^{(p)}_i$ for $i\in\N$
 in  appropriate regimes  when  $p\to 1$ and time tends to $\infty$. In this direction, 
we denote the mean reproduction of $Z^{(p)}_{\varnothing}$ by 
$$m_1(p)=\E(\xi)+p(1+\beta),$$ and recall that the process
$$W^{(p)}_{\varnothing}(t) = \e^{-m_1(p)t}Z^{(p)}_{\varnothing}(t)\,, \qquad t\geq 0$$ is a martingale with terminal value denoted by  $W^{(p)}_{\varnothing}(\infty)$. 
For each fixed $t\geq 0$, we have $\lim_{p\to 1}W^{(p)}_{\varnothing}(t)= W(t)$, and on the other hand, we know  that $\lim_{t\to \infty}W(t)= W(\infty)$ in $L^2$. As a matter of fact, we have a stronger uniform convergence. 

\begin{lemma}\label{Le3} It holds that
$$\lim_{p\to 1, t\to \infty} \E_z\left(\sup_{s\geq t} |W^{(p)}_{\varnothing}(s)-W(\infty)|^2\right)=0\,.$$
\end{lemma}

\proof  Note that for $1/2\leq p < 1$, we have $m_1(p)\geq \frac{1}{2}m_1$ and the second moment of $\xi+\epsilon_p (1+\beta)$ is at most $m_2$. 
We deduce from Lemma \ref{L1} applied to the branching process $Z^{(p)}_{\varnothing}$  that for every fixed $\varepsilon>0$,
we can find $t_{\varepsilon}<\infty$ such that
\begin{equation}\label{Eq3}
\E_z\left(\sup_{s\geq t_{\varepsilon}} |W^{(p)}_{\varnothing}(s)-W^{(p)}_{\varnothing}(\infty)|^2\right) \leq \varepsilon\qquad \hbox{for all }p\in[1/2,1]\,.
\end{equation}

We next claim that
\begin{equation}\label{Eq4}
\lim_{p\to 1} \E_z(|W^{(p)}_{\varnothing}(t_{\varepsilon})-W(t_{\varepsilon})|^2)=0\,.
\end{equation}
Indeed, recall that ${b}^{(p)}_1$ denotes the first birth time of a mutant population. Plainly,   $\lim_{p\to 1}{b}^{(p)}_1=\infty$ in probability, and the probability of the event $\{t_{\varepsilon}\geq {b}^{(p)}_1\}$ can be made as small as we wish by choosing $p$ sufficiently close to $1$. 
On the one hand, as
$Z^{(p)}_{\varnothing}(t_{\varepsilon})\leq Z(t_{\varepsilon})$, we have
$$\E_z(|W^{(p)}_{\varnothing}(t_{\varepsilon})-W(t_{\varepsilon})|^2, t_{\varepsilon}\geq {b}^{(p)}_1 )
\leq (\e^{2(m_1-m_1(p))t_{\varepsilon}}+1)\E_z(|W(t_{\varepsilon})|^2, t_{\varepsilon}\geq {b}^{(p)}_1 )\,,$$
and the right-hand side goes to $0$ as $p\to 1$. 
On the other hand, on the event $\{t_{\varepsilon}<  {b}^{(p)}_1\}$, we have
$Z^{(p)}_{\varnothing}(t_{\varepsilon})=Z(t_{\varepsilon})$ and hence
$W^{(p)}_{\varnothing}(t_{\varepsilon})=\e^{(m_1-m_1(p))t_{\varepsilon}}W(t_{\varepsilon})$. This yields 
$$\E_z(|W^{(p)}_{\varnothing}(t_{\varepsilon})-W(t_{\varepsilon})|^2, t_{\varepsilon}<  {b}^{(p)}_1 )
\leq (\e^{(m_1-m_1(p))t_{\varepsilon}}-1)^2\E_z(|W(t_{\varepsilon})|^2)\,,$$
and again the right-hand side goes to $0$ as $p\to 1$. This establishes \eqref{Eq4}. 

Combining \eqref{Eq3} and \eqref{Eq4}, we get 
$$\limsup_{p\to 1}\E_z(|W(\infty)-W^{(p)}_{\varnothing}(\infty)|^2)\leq 4\varepsilon\,,$$
and since $\varepsilon>0$ can be chosen arbitrarily small, we have in fact
\begin{equation}\label{Eq5}
\lim_{p\to 1}\E_z(|W(\infty)-W^{(p)}_{\varnothing}(\infty)|^2)=0\,.
\end{equation}
Plugging this in \eqref{Eq3}, we conclude that 
$$
\limsup_{p\to 1} \E_z\left(\sup_{s\geq t_{\varepsilon}} |W^{(p)}_{\varnothing}(s)-W(\infty)|^2\right) \leq \varepsilon\,,
$$
which is equivalent to our statement. \QED

We next turn our attention to the asymptotic behavior of the birth times ${b}^{(p)}_i$ for $i=1,2, \ldots$ of the different types with a single mutation.

  \begin{lemma}\label{Le4} As $p\to 1$, the sequence
   $$ \frac{1-p}{m_1} W(\infty)\,  \exp\left( m_1(p){b}^{(p)}_i\right)\,,\qquad i\geq 1$$
  converges in the sense of finite-dimensional distributions towards
$$S_i:={\bf e}_1+\cdots + {\bf e}_i\,,\qquad i\geq 1\,,$$
where $({\bf e}_i)_{i\in\N}$ denotes a sequence of i.i.d. standard exponential variables.
\end{lemma} 
We stress that, thanks to Lemma \ref{Le2}, the sequence above is independent of the processes 
 $(Z^{(p)}_i({b}^{(p)}_i+t): t\geq 0)$ for $i\geq 1$. This observation will be important later on. 
 
\proof
Define  
$$I^{(p)}(t):=(1-p)\int_{0}^t Z^{(p)}_{\varnothing}(s)\d s\,, \qquad t\geq 0\,.$$ The random map $I^{(p)}: [0,\infty)\to [0,\infty)$ is a.s. bijective, and we denote its inverse by $J^{(p)}$. 
It follows immediately from the description of the population system that 
if we time-change the  process $t \mapsto M^{(p)}(t)$ which counts the number of types with a single mutation, by 
$J^{(p)}$, then we get another counting process $t\mapsto M^{(p)}\circ J^{(p)}(t)$ with unit jump rate.
In other words, $M^{(p)}\circ J^{(p)}$ is a standard Poisson process, and therefore the sequence of its jump-times is given by a random walk $S$ with exponentially 
distributed steps with unit mean\footnote{
This random walk depends on the parameter $p$, however since only its law is relevant in this proof, this will be omitted from the notation for simplicity.}.
Since the birth-times ${b}^{(p)}_i$ for $i\geq 1$ are the jump-times of $M^{(p)}$, this yields
$$I^{(p)}({b}^{(p)}_i)=S_i\,,\qquad i\geq 1\,.$$

We now only need to estimate $I^{(p)}(t)$ as both $p\to 1$ and $t\to \infty$. 
In this direction observe from the triangle inequality that
\begin{eqnarray*}
 \left | I^{(p)}(t) - \frac{1-p}{m_1(p)}(\e^{m_1(p) t}-1)W^{(p)}_{\varnothing}(\infty) \right | 
 &\leq& (1-p) \int_{0}^t | Z^{(p)}_{\varnothing}(s)-W^{(p)}_{\varnothing}(\infty)\e^{m_1(p) s} |\d s \\
 &=& (1-p)\int_{0}^t  | W^{(p)}_{\varnothing}(s)-W^{(p)}_{\varnothing}(\infty)| \e^{m_1(p) s} \d s\\
 &\leq& (1-p)A^{(p)} \e^{m_1(p)2t/3}\,,
 \end{eqnarray*}
 where
$$A^{(p)}:=\frac{3}{2m_1(p)}\sup_{s\geq 0} \e^{m_1(p) s/3}| W^{(p)}_{\varnothing}(s)-W^{(p)}_{\varnothing}(\infty)|\,. $$
Recall from Lemma \ref{L1} that the variables $A^{(p)}$
 are bounded in $L^2(\P)$
for $1/2\leq p < 1$; we deduce that
$$\lim_{t\to \infty} \E\left( \sup_{s\geq t} \left| \frac{\e^{-m_1(p) s}}{1-p}
I^{(p)}(s) - \frac{W^{(p)}_{\varnothing}(\infty)}{m_1(p)} \right| ^2\right) =0\qquad \hbox{uniformly in  $1/2\leq p < 1$.}$$
Recall also from \eqref{Eq5} that $\lim_{p\to 1-}W^{(p)}_{\varnothing}(\infty)=W(\infty)$ 
in $L^2(\P)$, where $W(\infty)$  is strictly positive a.s., and note that $b^{(p)}_i\to \infty$ in probability for every $i\geq 1$. 
It follows  that 
$$S_i=I^{(p)}({b}^{(p)}_i)\sim \frac{1-p}{m_1(p)}\, \e^{m_1(p) {b}^{(p)}_i}W(\infty)
\qquad \hbox{in probability, }$$
and clearly we may replace $m_1(p)$ by $m_1$  in the fraction above. \QED

We have all the technical ingredients to establish the main result of this section, but we still need some additional notation.
For each $p\in(0,1)$, consider a random time $\tau^{(p)}$ such that 
\begin{equation}\label{Eq6}
\lim_{p\to 1} \left( m_1(p) \tau^{(p)} + \ln (1-p)\right) =\infty\qquad \hbox{in probability.}
\end{equation}
Let  $W'(\infty)$ be a variable distributed as the terminal value of the martingale $W(t)=\e^{-m_1t}Z(t)$ where the starting point  is now  $Z(0)=1+\beta$. 
We introduce $(W'_i(\infty): i\geq 1)$ a sequence of i.i.d.  copies of $W'(\infty)$. 
We finally recall that $(S_k: k\geq 0)$ denotes a random walk with i.i.d. steps distributed according to the standard exponential law.
We implicitly assume that $(W'_i(\infty): i\geq 1)$ and  $(S_k: k\geq 0)$ are independent. 

\begin{theorem} \label{T2}
 As $p\to 1$, the sequence 
$$\left( \frac{ \e^{-m_1(p) \tau^{(p)}}}{(1-p)W(\infty)}  Z_i^{(p)}(\tau^{(p)}):i\geq 1\right)$$
converges in the sense of finite-dimensional distributions  towards
$$\left(\frac{W'_i(\infty)}{m_1 S_i} : i\geq 1 \right)\,.$$

\end{theorem}

\proof Recall that ${b}^{(p)}_i$ denotes the instant of the $i$-th mutation in  the sub-population with the ancestral type, and set for $i\geq 1$ and $t\geq 0$, 
$$W^{(p)}_i(t) = \e^{-m_1(p)t}Z^{(p)}_i(t+{b}^{(p)}_i)\,.$$ 

Fix $z>0$. By Lemma \ref{Le3}, for every continuous $f:\R\to\R$ bounded in absolute value by $1$ and every $\varepsilon>0$, there exist $t(f,\varepsilon)$ and $p(f,\varepsilon)$ such that if $p(f,\varepsilon)<p\leq 1$ then$$
\E_z\left( \sup_{ t_1,t_2\geq t(f,\varepsilon) }\left| f(W^{(p)}_{\varnothing}(t_1))-f(W^{(p)}_{\varnothing}(t_2)) \right| \right)\leq \varepsilon
$$and
$$
\left|\E_z\left( f(W^{(p)}_{\varnothing}(t(f,\varepsilon)))\right) -\E_z\left(f(W(\infty)) \right)\right|\leq \varepsilon
$$
Without loss of generality, we may also assume that the same inequalities hold with 
$W^{(p)}_{\varnothing}$  and $W(\infty)$ replaced by $W^{(p)}_i$ and $W_i'(\infty)$ for any $i\in\N$, since this amounts to taking $z=1+\beta$. 

Consider then for each $i\geq 1$, a family of random times $(t^{(p)}_i)_{0<p<1}$, such that
$\lim_{p\to 1}t^{(p)}_i = \infty$ in probability.  Since we can guarantee that $\P_z\left(t^{(p)}_i\leq t(f,\varepsilon)\right)\leq \varepsilon$ for $p\geq p(i,f,\varepsilon)$, we see that
$$
\left|\E_z\left( f(W^{(p)}_i(t^{(p)}_i))\right) -\E\left(f(W'_i(\infty)) \right)\right|\leq 2\varepsilon
$$if $p\geq p(i,f,\varepsilon)$. The independence of the $W^{(p)}_i$ is seen from Lemma \ref{Le2}, and we have deduced the weak convergence (in the sense of finite dimensional distributions)
%
%
%
%
%
%
%
$$\left( \e^{-m_1(p)t^{(p)}_i}Z^{(p)}_i(t^{(p)}_i+{b}^{(p)}_i): i\geq 1 \right)
 \ \Longrightarrow \ (W'_i(\infty) : i\geq 1)\,.$$
 Recall from Lemma \ref{Le4} that we have also
 $$\left(  \frac{1}{(1-p)W(\infty)}\,  \exp\left( -m_1(p){b}^{(p)}_i\right):i\geq 1\right)
  \ \Longrightarrow \ \left(\frac{1}{m_1 S_i}: i\geq 1\right)\,.$$ 
  More precisely, we deduce from Lemma \ref{Le2} that these two weak convergences hold jointly,
  provided that we take the sequences $(W'_i(\infty) : i\geq 1)$ and $(S_i: i\geq 1)$ to be  independent. 
   
 To complete the proof, it now suffices to set  $t^{(p)}_i=\tau^{(p)}-{b}^{(p)}_i$ for $i\geq 1$ and take the product of  the preceding weak limits. Note that Lemma \ref{Le4}
 and the assumption \eqref{Eq6} ensure that indeed $\lim_{p\to 1}t^{(p)}_i = \infty$. \QED


\end{section}

\begin{section}{Combining preferential attachment with percolation}
\label{percolationSection}
We are now able to start investigating the question which has motivated this article.
It  is convenient for our purpose to work with a continuous version of the preferential attachment algorithm, in the sense that we shall grow a scale-free tree in continuous time. That is, we start at time $0$ from the tree on $\{0,1\}$, and 
once the random tree with size $n\geq 2$ has been constructed, we equip each vertex $i\in\{0,\ldots,n\}$ with an exponential clock $\zeta_i$ with parameter $d_n(i)+\beta$, where $d_n(i)$ denotes the current degree of $i$, independently of the other vertices. Then  the next vertex $n+1$ is attached after  time $\min_{i\in\{0,\ldots,n\}}\zeta_i$ at the vertex
$v_n=\arg\min_{i\in\{0,\ldots,n\}} \zeta_i$. 
Recall that the sum of the degrees of a tree with $n+1$ vertices is $2n$, so $\min_{i\in\{0,\ldots,n\}}\zeta_i$
is exponentially distributed with parameter $2n+ \beta(n+1)$.

Denote by $T(t)$ the tree which has been constructed at time $t$, and by $|T(t)|$ its size, i.e. its number of vertices.
It should be plain that if we define 
$$\tau_n=\inf\{t\geq 0: |T(t)|=n+1\}\,,$$ then $T(\tau_n)$ is a version of a scale-free tree of size $n+1$, $T_n$.
The process of the size $|T(t)|$ of $T(t)$ is clearly Markovian; however it will more convenient in practice to work with a linear transformation of it, namely
$$Y(t)= 2(|T(t)|-1)+\beta|T(t)|\,,æ\qquad t\geq 0.$$
In particular $Y(0)=2+2\beta$.

\begin{lemma}\label{Le5} The process  $Y$ is a pure birth branching process, that  has only jumps of size $2+\beta$, and with unit birth rate per unit population size. Equivalently, 
$(2+\beta)^{-1}Y$ is a Yule branching process in continuous space with birth rate $2+\beta$ per unit population size. 
\end{lemma} 

\proof The sum of degrees of vertices in $T(t)$ is
$2(|T(t)|-1)$. Because when this tree  has size $n$, the next vertex $n+1$ is incorporated at rate $2(n-1)+ \beta n$, which yields an increase of  $Y$ by $2+\beta$,
we see that $Y$   is a branching process in continuous space and time with unit rate of birth per unit population size, and reproduction law  given by the Dirac point mass at $2+\beta$. Normalizing $Y$ by a factor $(2+\beta)^{-1}$, we recognize the dynamics of a Yule process. 
 \QED

We next superpose Bernoulli bond percolation to this construction by marking each edge $e_j$ connecting a vertex $j\geq 1$ to its parent $v_j$ with  an independent uniform random variable $U_j$. The parameter $p\in(0,1)$ being fixed, we may imagine that  $e_j$  is cut at its midpoint when the mark $U_j>p$ and remains intact otherwise.
We write $T^{(p)}(t)$ for the resulting combinatorial structure at time $t$. That is  $T^{(p)}(t)$ has the same set of vertices as $T(t)$, its set of intact edges is the subset of the edges $e_j$ of $T(t)$ such that $U_j\leq p$, and further $T^{(p)}(t)$ may have half-edges which should be viewed as stubs attached to some vertices and correspond to edges of $T(t)$ which have been cut in two.
The point in cutting  rather than removing edges is that the former procedure preserves the degrees of vertices, where the degree of a vertex is defined as the sum of the intact edges and half-edges attached to it.

The percolation clusters of $T^{(p)}(t)$ are the subtrees of $T(t)$ formed by the subsets of vertices which can be connected by a path of intact edges. We write $T_0^{(p)}(t), T_1^{(p)}(t), \ldots$
for the sequence of subtrees at time $t$, where the enumeration follows the increasing order of their birth times, and with the convention that $T_j^{(p)}(t)=\emptyset$ when the number of edges that have been cut at time $t$ is less than $j$. Specifically, if $j$ is the label of the $i$-th  variable $U_j$ to be greater than $p$, then $T^{(p)}_i(t)$ is the combinatorial structure spanned by the vertices that can be joined by a path of intact edges to the vertex $j$. In particular $T_0^{(p)}(t)$ denotes the subtree at time $t$ which contains the vertex $0$; it shall play a special role in our analysis. 

We write $H^{(p)}_i(t)$ for the number of half-edges pertaining to the $i$-th subtree at time $t$, so that  $2(|T_i^{(p)}(t)|-1) +  H^{(p)}_i(t)$ is the sum of the degrees of vertices of the $i$-th subtree. 
We stress that
$$\sum_{i\geq 0} |T^{(p)}_i(t)| = |T(t)| \quad\hbox{and}\quad \sum_{i\geq 0} (2(|T^{(p)}_i(t)|-1)+ H^{(p)}_i(t)) = 2(|T(t)|-1)\,.$$
If we set
$$ Y^{(p)}_i(t)=2(|T_i^{(p)}(t)|-1) +  H^{(p)}_i(t)+ \beta |T_i^{(p)}(t)|\,,\qquad t\geq 0\,,$$
then we see from above that 
\begin{equation}\label{Eq7}
\sum_{i\geq 0} Y^{(p)}_i(t) = Y(t)\,.
\end{equation}

The connexion with the system of branching processes with neutral mutations of the preceding section should be clear. Specifically, imagine that at some given time $t$, the state of the the process ${\bf Y}^{(p)}= (Y^{(p)}_j : j\geq 0)$ is given by $(y_0, y_1, \ldots)$, and write $y=y_0+y_1+\cdots$. In particular the current size of the growing tree is $|T(t)|=(y+2)/(2+\beta)$ and we know from Lemma \ref{Le5} that the next vertex will be incorporated after an exponential time with parameter $y$.
The probability that the edge corresponding to this new vertex has its other extremity in the $i$-th subtree $T^{(p)}_i(t)$ is 
$$\frac{2(|T_i^{(p)}(t)|-1) +  H^{(p)}_i(t)+ \beta |T_i^{(p)}(t)|}{y}=\frac{y_i}{y}\,,$$
independently of the waiting time. Finally, the probability that this edge is intact is $p$, independently of the preceding variables. We thus see from basic properties of independent exponential variables that ${\bf Y}^{(p)}$ has  the same random evolution as the system ${\bf Z}^{(p)}$ of branching processes with neutral mutations of Section  \ref{branchingSection} when the reproduction law $\nu$ is simply given by the Dirac mass at $2+\beta$. 
 In this setting, $Y^{(p)}_{0}$ corresponds to $Z^{(p)}_{\varnothing}$, the sub-population with the ancestral type of Section 2. There is however an important difference between the way clusters and sub-populations are labeled that should be stressed to avoid a possible confusion. More precisely, 
 the families ${\bf Z}^{(p)}= (Z^{(p)}_u: u\in {\mathbf U})$ and ${\bf Y}^{(p)}= (Y^{(p)}_j : j\geq 0)$ do represent the same process, 
and in particular  $(Z^{(p)}_i: i\in\N)$ is only a sub-sequence of $(Y^{(p)}_j : j\in\N)$ which corresponds to subtrees at distance $1$ from the root-subtree $T^{(p)}_0$.  Recall that focussing on sub-populations with a single mutation is crucial to ensure  the validity of Lemma \ref{Le2}.

Define the generation of a vertex  as the number of edges $e$ on the branch from this vertex to the root $0$ which have a mark $U_e>p$ (in other words, this is the number of cuts on that branch). In particular vertices of $T^{(p)}_0$ have generation $0$, and those of $T^{(p)}_1$ have generation $1$.  We then set $\rho(i)=j$ for $i\geq 1$ when the $j$-th subtree of $T^{(p)}(t)$ is the $i$-th subtree of the first generation, where, as usual,  subtrees in a family are enumerated  according to the increasing order of their birth times. In particular, we always have $\rho(1)=1$ and the sequence $(\rho(i): i\geq 1)$ is strictly increasing. 
The following claim should be plain from the discussion above.

\begin{corollary}\label{Co1} In the notation of Section 2, take  $\xi\equiv 1$ and $z=2+2\beta$. Then the families 
$$(Y, Y^{(p)}_{0}, Y^{(p)}_{\rho(1)}, Y^{(p)}_{\rho(2)}, \ldots) \quad \hbox{and}\quad (Z, Z^{(p)}_{\varnothing}, Z^{(p)}_{1}, Z^{(p)}_{2}, \ldots)$$
have the same distribution. 
 \end{corollary} 
 In the sequel, it will be convenient to agree that the two families in the statement above are actually the same (not merely are identical in law).
Recall also that the algorithm with preferential attachment is run until time
$$\tau_n=\inf\{t\geq 0: |T(t)|=n+1\}= \inf\{t\geq 0: Y(t)=2n+\beta (n+1)\}$$ when the structure has size $n+1$. 

{\bf We henceforth assume that the percolation parameter $p=p(n)$ fulfills \eqref{Eq1}.} The motivation for this choice  stems from the next statement, which shows that both the root-cluster and its complement are then macroscopic (i.e. of size of order $n$). 
For the sake of simplicity, we shall frequently write $p$ rather than $p(n)$, omitting the integer $n$ from the notation. Recall that $\alpha= (1+\beta)/(2+\beta)$.

\begin{corollary} \label{Co2}
We have
$$\lim_{n\to \infty} \frac{Y^{(p)}_{0}(\tau_n)}{n} =(2+\beta)\e^{-\alpha c }$$
in probability.
\end{corollary} 
\proof   We know from Corollary \ref{Co1} and Lemmas \ref{L1} and  \ref{Le3} that 
$$\lim_{n\to \infty} \e^{-m_1\tau_n} Y(\tau_n)  = \lim_{n\to \infty} \e^{-m_1(p)\tau_n} Y^{(p)}_{0}(\tau_n) = W(\infty)\qquad \hbox{in probability,}$$
with $m_1=2+\beta$ and $m_1(p)=1+p(1+\beta)$. 
By the definition of $\tau_n$, we have $Y(\tau_n)=2n+\beta (n+1)$,
 hence
$$\e^{-(2+\beta)\tau_n}\sim \frac{W(\infty)}{(2+\beta)n}$$
and, {\it a fortiori}, $\tau_n\sim (2+\beta)^{-1} \ln n$. Our claim follows since
$$m_1-m_1(p)=(1-p(n))(1+\beta)\sim \frac{(1+\beta)c}{\ln n}\,,$$
 thank to \eqref{Eq1}. \QED 

Next, let 
$$N^{(p)}(t)=\max\{j: T_j^{(p)}(t)\neq \varnothing\}$$ 
denote the number of subtrees at time $t$, discounting the root-subtree $T^{(p)}_0(t)$ containing $0$. Recall also that $M^{(p)}(t)$ denotes the number of sub-populations with a single mutation at time $t$, that is of subtrees at time $t$ which are at unit distance from  $T^{(p)}_0(t)$. We shall now observe that when $p$ is close to $1$, these two quantities coincide with high probability as long as they  are not too large. In this direction, recall  that $p=p(n)$ fulfills \eqref{Eq1}, and observe from Lemma \ref{Le4} that for each fixed $i\geq 1$,  the time ${b}^{(p)}_i$ of the $i$-th mutation within the sub-population with the ancestral type (i.e. the first instant when $M^{(p)}$ reaches $i$) fulfills 
$${b}^{(p)}_i= \frac{1}{m_1(p)}\ln \frac{1}{1-p(n)} + O(1) = \frac{\ln\ln n}{2+\beta}+O(1)\qquad \hbox{as }n\to \infty\,.$$

\begin{lemma} \label{Le6} Set $\Delta^{(p)}(t)=N^{(p)}(t)-M^{(p)}(t)$ for the number of subtrees at distance strictly greater than $1$ from the root-cluster at time $t$.
For every $r>0$, we have
$$\lim_{n\to \infty} \E\left(\Delta^{(p)}((2+\beta)^{-1} \ln \ln n + r)\right) =0\,.$$

\end{lemma} 

\proof Roughly speaking, the dynamics of ${\bf Y}^{(p)}$ show that  the counting process $N^{(p)}$ grows at rate $(1-p)Y$, which means rigorously that 
the predictable compensator of $N^{(p)}$ is absolutely continuous with density $(1-p)Y$. In other words,  $N^{(p)}(t)-(1-p)\int_{0}^tY(s)\d s$ is a martingale, and thus
$$ \E\left(N^{(p)}((2+\beta)^{-1} \ln \ln n + r) \right)  = (1-p) \int_{0}^{ (2+\beta)^{-1} \ln \ln n + r} \E(Y(s))\d s\,.$$

Similarly, the counting process $M^{(p)}$ grows at rate $(1-p)Y^{(p)}_{0}$, and 
$$ \E\left(M^{(p)}((2+\beta)^{-1} \ln \ln n + r) \right)  = (1-p) \int_{0}^{ (2+\beta)^{-1} \ln \ln n + r} \E(Y^{(p)}_{0}(s))\d s\,.$$

We deduce from Lemma \ref{Le5} that 
$$\E(Y(s))= 2(1+\beta)\e^{(2+\beta)s}\quad \hbox{and} \quad \E(Y^{(p)}_{0}(s))= 2(1+\beta)\e^{(2+\beta-(1-p)(1+\beta))s}\,,$$
 and our claim then follows from \eqref{Eq1}. \QED

Lemma \ref{Le6} entails in particular that for each fixed $k\geq 1$, the probability that the $k$-tuple of processes 
$(Y^{(p)}_i)_{1\leq i \leq k}$ and $(Z^{(p)}_i)_{1\leq i \leq k}$ coincide tends to $1$ as $n\to \infty$. 
This enables us to deduce the asymptotic behavior of the former from Theorem \ref{T2}. We shall use the same notation as there, specialized to the setting of this present section. 
That is, $W'(\infty)$ denotes the terminal value of the martingale $\e^{-(2+\beta)t}Y(t)$ given $Y(0)=1+\beta$, 
 $(W'_i(\infty))_{i\geq 1}$  is  a sequence of i.i.d. copies of $W'(\infty)$, 
and $(S_i)_{i\geq 1}$ an independent random walk whose steps have the standard exponential distribution.

\begin{corollary} \label{Co3} 
The sequence 
$$\left( \frac{  \ln n }{n}  Y^{(p)}_i(\tau_n):i\geq 1\right)$$
converges in the sense of finite-dimensional distributions as $n\to \infty$  towards
$$\left(c \e^{-\alpha c }\frac{W'_i(\infty)}{S_i} : i\geq 1 \right)\,.$$
\end{corollary}

 \proof Recall that  $(m_1-m_1(p))\tau_n\to \alpha c $, as proved in Corollary \ref{Co2}. Hence
 $$\exp(-m_1(p)\tau_n)=  \exp((m_1-m_1(p))\tau_n)\exp(-m_1\tau_n) \sim  \e^{\alpha c } \frac{W(\infty)}{(2+\beta)n}\,,$$
and our claim follows from Theorem \ref{T2} specified  in the present setting with $\tau^{(p)}=\tau^{(p(n))}=\tau_n$. \QED

Next, we easily translate the above limit theorem for the branching processes $Y^{(p)}_i$ in terms of the sizes of the subtrees  listed in the increasing order of their ages. 

  \begin{corollary} \label{Co4} We have
  $$\lim_{n\to \infty} n^{-1} |T^{(p)}_{0}(\tau_n)| = \e^{-\alpha c }$$
and the sequence 
$$\left( \frac{\ln n  }{n}  |T^{(p)}_{i}(\tau_n)| :i\geq 1\right)$$
converges as $n\to \infty$,  in the sense of finite-dimensional distributions,  towards
$$\left(c \e^{-\alpha c }  \frac{W'_i(\infty)}{(2+\beta)S_i} : i\geq 1 \right)\,.$$
\end{corollary}

\proof We  focus on the second claim, the proof of the first being similar (and easier) using Corollary \ref{Co2} in place of Corollary \ref{Co3}. 

From Corollary \ref{Co3}, it suffices to show that
$$Y^{(p)}_i(\tau_n) \sim (2+\beta)  |T^{(p)}_i(\tau_n) |\,,$$
and for this, that the number of half-edges pertaining to the $i$-th sub-tree  fulfills 
\begin{equation}\label{Es}
H^{(p)}_i(\tau_n) = o\left(Y^{(p)}_i(\tau_n)\right)\,.
\end{equation}

In this direction, recall that  the $i$-th  jump time $\gamma^{(p)}_i:=\inf\{t\geq 0: N^{(p)}(t)=i\}$ of the process  $N^{(p)}$ that counts the number of subtrees as time passes, is a stopping time which corresponds to the birth-time of the $i$-th subtree $T^{(p)}_i$. We observe from the dynamics described at the beginning of this section and the strong Markov property  that the process
$$H^{(p)}_i(\gamma^{(p)}_i+t) - (1-p)\int_{0}^t Y^{(p)}_i(\gamma^{(p)}_i+s) \d s\,, \qquad t\geq 0 $$
is a martingale. 
Similarly, 
$$Y^{(p)}_i(\gamma^{(p)}_i+t) - (1-p + p(2+\beta)) \int_{0}^t Y^{(p)}_i(\gamma^{(p)}_i+s) \d s\,, \qquad t\geq 0 $$
is also a  martingale.  It follows that  
$$L^{(p)}(t):=H^{(p)}_i(\gamma^{(p)}_i+t)  - \frac{1-p}{1+p+p\beta} Y^{(p)}_i(\gamma^{(p)}_i+t)$$
is a martingale; note also that its jumps $|L^{(p)}(t)-L^{(p)}(t-)|$ have size at most $2+\beta$, independently of $p$. Since there are at most $n$ jumps up to time $\tau_n-\gamma^{(p)}_i$, 
the bracket of $L^{(p)}$ can be bounded by
$$[L^{(p)}]_{\tau_n-\gamma^{(p)}_i} \leq (2+\beta)^2n\,.$$
Hence
$$\E(|L^{(p)}(\tau_n-\gamma^{(p)}_i) - L^{(p)}(0)|^2)\leq  (2+\beta)^2n\,,$$
and in particular 
$$\lim_{n\to \infty}\E\left(\left |
\frac{ \ln n}{n} L^{(p)}(\tau_n-\gamma^{(p)}_i)\right|^2\right) =0\,.$$ 
The estimate 
 \eqref{Es} now follows readily from Corollary \ref{Co3} and the fact that $1-p(n)=o(1)$. 
\QED

Our final task is to deduce from Corollary \ref{Co3} a limit theorem for the sizes of the percolation clusters listed in the decreasing order of their sizes, rather than their ages. 
Roughly speaking, we shall check that the largest clusters are given by the older subtrees, in the sense that for every fixed $k$, with high probability when $n\to \infty$ and $\ell\to \infty$, 
the $k$ largest percolation clusters of $T^{(p)}(\tau_n)$ are to be found amongst the $\ell$ oldest subtrees $(T^{(p)}_i(\tau_n))_{0\leq i \leq \ell}$. 

Recall from Lemma \ref{Le4} that for each fixed $i\geq 1$,  the $i$-th oldest subtree $T^{(p)}_i(\tau_n)$ was born at time $(2+\beta)^{-1} \ln \ln n+O(1)$, and  from Corollary \ref{Co3} that its size is of order $n/\ln n$. We thus have to check that it is unlikely  to have at time  $\tau_n$ a subtree of size $\approx n/\ln n$ or greater, and which was born at a much later time than $(2+\beta)^{-1} \ln \ln n$.
Here is a formal statement, which is expressed for conveniency in terms of the processes $Y^{(p)}_k$.

\begin{lemma} \label{Le7} For every  $\varepsilon >0$, we have
$$\lim_{r\to \infty} \limsup_{n\to \infty} \P\left(\exists k\geq 1: Y^{(p)}_k((2+\beta)^{-1} \ln \ln n + r)=0 \hbox{ and } Y^{(p)}_k(\tau_n)> \varepsilon n/\ln n\right) =0\,.$$

\end{lemma} 

\proof Let $({\mathcal F}_t)_{t\geq 0}$ denote the natural filtration generated by the (continuous time version of) the algorithm with preferential attachment, including the uniform marks on the edges. The counting process $N^{(p)}$ is $({\mathcal F}_t)$-adapted, and its jump times $\gamma^{(p)}_k:=\inf\{t\geq 0: N^{(p)}(t)=k\}$ are stopping times that correspond to the birth-times of subtrees.

An application of the strong Markov property to the algorithm (recall also Lemma \ref{Le5}) shows that for each $k\geq 1$, the process $(2+\beta)^{-1}Y^{(p)}_k(\cdot + \gamma^{(p)}_k)$
is a Yule process with birth rate $2+\beta$ per unit population size, started from $(1+\beta)/(2+\beta)=\alpha\leq 1$, and independent of ${\mathcal F}_{\gamma^{(p)}_k}$. 
Plainly, the latter can be bounded from above by a Yule process with the same birth rate and started at $1$, in particular 
$(2+\beta)^{-1}Y^{(p)}_k(u + \gamma^{(p)}_k)$ is stochastically bounded from above by a geometric variable with parameter $\exp(-(2+\beta)u)$ (see, e.g. \cite{MR2047480} on page 109). That is, the tail distribution of $(2+\beta)^{-1}Y^{(p)}_k(u + \gamma^{(p)}_k)$ admits the bounds $\ell\mapsto (1-\exp(-(2+\beta)u))^{\ell+1}$.

 It is convenient to write $r_n=(2+\beta)^{-1} \ln \ln n + r$ and $s_n=(2+\beta)^{-1}  \ln n + s$.
Fix $s\geq 0$ arbitrary large, and consider the number of processes $Y^{(p)}_{\cdot}$ which are born after time $r_n$ and reach a size greater than $\varepsilon n/\ln n$
at time $s_n$, namely
$$X_n=\sum_{k=1}^{\infty} {\bf 1}_{\{{r_n\leq \gamma^{(p)}_k \leq s_n }\}} {\bf 1}_{\{Y^{(p)}_k(s_n)>\varepsilon n/\ln n\}}= \int_{r_n}^{s_n}  {\bf 1}_{\{Y^{(p)}_{N^{(p)}(t)}(s_n)>\varepsilon n/\ln n\}} \d N^{(p)}(t)
.$$
The preceding observations entail that
\begin{eqnarray*}
\E(X_n) &\leq& \E\left(\int_{r_n}^{s_n }  (1-\exp(-(2+\beta)(s_n-t)))^{(2+\beta)^{-1}\varepsilon n/\ln n} \d N^{(p)}(t)\right) \\
&\leq& \E\left(\int_{r_n}^{s_n }  \exp\left[-\frac{\varepsilon n}{(2+\beta) \ln n}\exp(-(2+\beta)(s_n-t))\right] \d N^{(p)}(t)\right)\,,
\end{eqnarray*}
where in the second line, we used the inequality $(1-x)^a \leq \exp(-ax)$. 

Next recall that the counting process $N^{(p)}$ has a predictable compensator which is absolutely continuous with density $(1-p)Y$. This enables us to express the last quantity above as 
$$(1-p)\E\left(\int_{r_n}^{s_n } Y(t) \exp\left[-\frac{\varepsilon n}{(2+\beta)\ln n}\exp(-(2+\beta)(s_n-t))\right]\, \d t\right)\,.$$
Recall also that
$\E(Y(t))= 2(1+\beta)\e^{(2+\beta)t}$; we arrive at
\begin{eqnarray*}
\E(X_n) &\leq& 2(1+\beta) (1-p)\int_{r_n}^{s_n } \e^{(2+\beta)t} \exp\left[-\frac{\varepsilon n}{(2+\beta)\ln n}\exp(-(2+\beta)(s_n-t))\right]\, \d t \\
&=& \frac{2(1+\beta)}{2+\beta} (1-p)\int_{\e^{(2+\beta)r_n} }^{\e^{(2+\beta)s_n}  }  \exp\left[-x\frac{\varepsilon n}{(2+\beta)\ln n}\exp(-(2+\beta)s_n) \right]\, \d x \\
&\leq &2(1+\beta) (1-p) \frac{\ln n}{\varepsilon n}\exp((2+\beta)s_n)\exp\left[-\e^{(2+\beta)r_n} \frac{\varepsilon n}{(2+\beta)\ln n}\exp(-(2+\beta)s_n)  \right]\,.
\end{eqnarray*}

Plugging into the last expression the values of $r_n$ and $s_n$ and using \eqref{Eq1}, we conclude that 
$$\limsup_{n\to \infty}\E(X_n) \leq c' \varepsilon^{-1}\e^{(2+\beta)s} \exp(-\varepsilon(r-s)/(2+\beta))\,,$$
where $c'$ is some constant depending only on $\alpha,\beta$ and $c$. This quantity goes to $0$ as $r\to \infty$, for every fixed $\varepsilon$ and $s$, and this entails in particular that 
$$\lim_{r\to \infty} \limsup_{n\to \infty} \P\left(\exists k: Y^{(p)}_k(r_n)=0 \hbox{ and } Y^{(p)}_k(s_n)> \varepsilon n/\ln n\right) =0\,.$$
Since (from Lemma \ref{Le5})
$$\lim_{s\to \infty}\P(\tau_n> s_n=(2+\beta)^{-1}  \ln n + s)=0\,,$$
this completes the proof. \QED

We are now able to complete this paper and establish Theorem \ref{T1}. Recall that the scale-free tree $T_n$ with $n+1$ vertices can be obtained as $T(\tau_n)$, and that
the first claim of  Theorem \ref{T1} is also the first claim of Corollary \ref{Co4}. The key issue for the second claim is the following. Corollary \ref{Co3} provides a limit theorem in the sense of finite-dimensional distributions for the normalized sequence of the sizes of the subtrees ordered by age, whereas Theorem \ref{T1} concerns ordering by size.

Lemma \ref{Le7} and the fact that $|T^{(p)}_k(t)|-1\leq (2+\beta)^{-1}Y^{(p)}_k(t)$ enable us to assert that if we use the notation $(x_i)^{\downarrow}$ for the decreasing rearrangement of a sequence of nonnegative real numbers $(x_i)$ which converges to $0$, then 
$$\left( \frac{\ln n  }{n}  |T^{(p)}_{i}(\tau_n)| :i\geq 1\right)^{\downarrow}$$
converges as $n\to \infty$,  in the sense of finite-dimensional distributions,  towards
$$\left(c \e^{-\alpha c }(2+\beta)^{-1}  \frac{W'_i(\infty)}{S_i} : i\geq 1 \right)^{\downarrow}\,.$$

So all what is needed now is to identify explicitly the distribution of the limit above. In this direction, we start specifying the law of the i.i.d. variables $W'_i(\infty)$. 
From Lemma \ref{Le5} and standard properties of Yule processes (recall also  the notation $\alpha=(1+\beta)/(2+\beta)$), we get that  that if $Y'$ denotes a version of the branching process $Y$ started from $Y'(0)=1+\beta$, then  
$$ \lim_{t\to \infty}\e^{-(2+\beta)t}Y'(t)= W'(\infty)  \quad\qquad \hbox{a.s. and in }L^2(\P)$$
where $W'(\infty)$ is a gamma variable with (shape and rate) parameter $(\alpha, 1/(2+\beta))$, that is
$$\P(W'(\infty)\in \d w)=   \Gamma(\alpha)^{-1} (2+\beta)^{-\alpha} w^{\alpha-1}\e^{-w/(2+\beta)} \d w\,,\qquad w>0\,.$$


Recall that $(S_i)_{i\geq 0}$ is a random walk with standard exponentially distributed steps, which is independent of the $(W'_i(\infty))_{i\geq 1}$. It follows that 
 $(S_i, W'_i)_{i\geq 1}$ can be viewed as the sequence of the atoms of a Poisson point process on $\R_+\times \R_+$ with intensity
$$\d s \otimes (2+\beta)^{-\alpha} \frac{w^{\alpha-1}}{\Gamma(\alpha)} \e^{-w/(2+\beta)} \d w\,.$$
The image of this measure by the map 
$$(s,w)\mapsto x=c \e^{-\alpha c }(2+\beta)^{-1}  \frac{w}{s}$$
is again a Poisson random measure, now with intensity $c \e^{-\alpha c } x^{-2}\d x$, which establishes the second part of Theorem \ref{T1}. 
Finally, the alternative description in terms of the inverses of the atoms belongs to the folklore of Poisson random measures (simply note that the image of $c \e^{-\alpha c} x^{-2}\d x$ by the map $x\mapsto 1/x$ is  $c \e^{-\alpha c} \d x$). \end{section}

\bibliography{GenBib}
\bibliographystyle{plainnat}
\end{document}